\documentclass[12pt]{article}

\usepackage{amsmath,amssymb,amsthm,bbm}
\usepackage{color,framed,accents,tikz}

\normalsize

\newtheorem{theorem}{Theorem}[section]
\newtheorem{lemma}[theorem]{Lemma}
\newtheorem{remark}[theorem]{Remark}

\newtheorem{corollary}[theorem]{Corollary}

\numberwithin{equation}{section}

\def\({\bigl(}   \def\){\bigr)}

\def\musize{\sigma}
\def\Hrda{{\vec{H}(\boldsymbol{d}^+,\boldsymbol{d}^-,\musize)}}

\def\calHdra{{\vec{\mathcal{H}}(\boldsymbol{d}^+,\boldsymbol{d}^-,\musize)}}

\def\dvec{{\boldsymbol{d}}}
\def\kvec{{\boldsymbol{k}}}

\def\kmax{k_{\max}}
\def\dmax{d_{\max}}
\def\smax{s_{\max}}
\def\tmax{t_{\max}}

\def\calB{\vec{\mathcal{B}}}  
\def\B{{\vec{B}}}  
\def\svec{\boldsymbol{s}}
\def\tvec{\boldsymbol{t}}
\def\gvec{\boldsymbol{g}}
\def\ellvec{\boldsymbol{\ell}}
\def\mplus{M^+}
\def\mminus{M^-}

\def\dfrac#1#2{\lower0.15ex\hbox{\large$\frac{#1}{#2}$}}
\def\E{{\mathbb{E}}}

\def\dfrac#1#2{\lower0.15ex\hbox{\large$\frac{#1}{#2}$}}
\let\originalleft\left
\let\originalright\right
\renewcommand{\left}{\mathopen{}\mathclose\bgroup\originalleft}
\renewcommand{\right}{\aftergroup\egroup\originalright}

\allowdisplaybreaks

\title{Enumeration of dihypergraphs with\\ specified degrees and edge types}

\author{Catherine Greenhill\thanks{Supported by Australian Research Council Discovery Project DP250101611.}\\
\small School of Mathematics and Statistics\\[-0.5ex]
\small UNSW Sydney\\[-0.5ex]
\small NSW 2052, Australia\\
\small \texttt{c.greenhill@unsw.edu.au}
\and
Tam{\' a}s Makai\\
\small Department of Mathematics\\[-0.5ex]
\small LMU Munich\\[-0.5ex]
\small Munich 80333, Germany\\
\small \texttt{makai@math.lmu.de}
}

\begin{document}

\maketitle

\begin{abstract}
A directed hypergraph (dihypergraph) consists of a set of vertices and a set of hyperarcs, where each hyperarc is partitioned into a head and a tail. Directed hypergraphs
are useful in many applications, including the study of chemical reactions 
or relational databases. We provide asymptotic formulae for the number
of directed hypergraphs with given in-degree sequence,  out-degree sequence, and with the
head and tail sizes of all hyperarcs specified. Our formulae hold when none of the following parameters are too large: the maximum out-degree, the maximum in-degree, the maximum head size and the maximum tail size. 
\end{abstract}

\section{Introduction}\label{s:intro}

Given a finite set of vertices $V$,
a \emph{hyperarc} is an ordered pair $e=(e^+,e^-)$ of nonempty disjoint subsets of $V$. We say that
$e^+$ is the \emph{tail} and $e^-$ is the \emph{head} of $e$.
A \emph{directed hypergraph}  (dihypergraph) $H=(V,E)$ consists of a finite set $V$ of vertices
and a set $E$ of hyperarcs.  It follows from this definition that a directed hypergraph does not contain any repeated hyperarcs or ``loops'', which are hyperarcs that contain a repeated vertex.
Directed hypergraphs arise in many applications, for example when
modelling chemical reactions~\cite{ER-V,GBSJR} or in the study of relational
databases and satisfiability formulae~\cite{AL}.
There is increasing interest in directed hypergraphs (or \emph{dihypergraphs}, for short) from the network
science community, for example~\cite{KS,MKKS}. 

Our dihypergraphs definition
follows that used by Gallo, Longo, Pallottino and Nguyen~\cite{GLPN} and Klamt, Haus and Theis~\cite{KHT}, except that~\cite{GLPN} allows the tail or head of a hyperarc to be empty. Our definition is more general than some in the literature. For example,
the survey paper of Ausiello and Laura~\cite{AL} uses a definition
which restricts $e^-$ to always consist of a single vertex, called the
head of $e$. Other works using this definition are~\cite{BJJ,FKK}. 
Gallo et al.~\cite{GLPN} refer to such hypergraphs as \emph{B-hypergraphs},
and define \emph{F-hypergraphs} to be 
dihypergraphs in which $|e^+|=1$ for every hyperarc $e$.
They also consider BF-hypergraphs in which every hyperarc satisfies either
$|e^-|=1$ or $|e^+|=1$. Thakur and Tripathi~\cite{TT} compare these three
dihypergraph models. 

Qian~\cite{qian} gave an exact enumeration result 
for unlabelled B-hypergraphs and unlabelled $k$-uniform B-hypergraphs
with a given number of vertices.  Here ``unlabelled'' means that these dihypergraphs
were counted up to isomorphism.  We are not aware or any other enumeration
result for dihypergraphs in the literature.

It is useful to have asymptotic enumeration formulae for families of
combinatorial structures.  For example, asymptotic enumeration results for graphs
and bipartite graphs
with given degrees have been applied both within mathematics~\cite{BLW,FKNSS} and in
other disciplines~\cite{blanchet,MSW}.
Asymptotic enumeration formulae for directed graphs with given degree sequence
have been obtained by Liebenau and Wormald~\cite{LW} for a wide range
of degrees, and by Greenhill and McKay~\cite{GM09} for dense degree 
sequences, see also Barvinok~\cite{barvinok}.

Our aim in this paper is to use the switching method to
provide an asymptotic formula for the number of
dihypergraphs when the degrees and hyperarc sizes are not too large. 
We will use the connection between dihypergraphs and directed bipartite graphs,
which allows us to make use of existing enumeration formulae for sparse
bipartite graphs.

\subsection{Notation and statement of our results}\label{s:statements}

Let $\mathbb{N} = \{0,1,2,\ldots\}$ denote the natural numbers
and let $\mathbb{Z}_+$ denote the positive integers.
Given a directed hypergraph $H = (V,E)$ with $|V|=n$, 
let $m:=|E|$ be the number
of hyperarcs in $H$.  Vertex $v\in V$ has out-degree 
$d^+_v$ and in-degree $d^-_v$ defined by
\[ d^+_v := |\{ e\in E\mid v \in e^+\}|,\qquad
   d^-_v := |\{ e\in E\mid v \in e^-\}|,\]
and the degree sequence of $H$ is $\dvec:=(\dvec^+,\dvec^-)$ where
\[ \dvec^+ := \big(d^+_1, d^+_2,\ldots, d^+_n\big),\qquad
   \dvec^- := \big(d^-_1, d^-_2,\ldots, d^-_n\big).\]
We also use a function $\musize$ to capture the number of vertices in the head and tail of each hyperarc. 
Specifically, $\musize:\mathbb{Z}_+^2\to\mathbb{N}$ is a function such that there are exactly $\musize(a^+,a^-)$ hyperarcs $e\in E$ with $|e^+|=a^+$ and $|e^-|=a^-$.  
Then the number of hyperarcs in the directed hypergraph is $m=\sum_{(a^+,a^-)\in \mathbb{Z}_+^2}\, \musize(a^+,a^-)$ and
\begin{equation}\label{eq:degedgesize}
 \mplus :=\sum_{v \in [n]} d_v^+ = \sum_{(a^+,a^-)\in\mathbb{Z}_+^2} a^+\, \musize(a^+,a^-), 
  \quad
    \mminus :=\sum_{v \in [n]} d_v^- = \sum_{(a^+,a^-)\in\mathbb{Z}_+^2} a^-\, \musize(a^+,a^-).
\end{equation}
Let $\calHdra$ be the set of directed hypergraphs with degrees $\dvec^+$, $\dvec^-$ and with head and tail sizes given by $\musize$, and let
$\Hrda:= |\calHdra|$. We want to approximate $\Hrda$ in the sparse
setting, where the maximum degrees and maximum head and tail sizes
are not too large.
To avoid trivialities we assume
throughout the paper that $\mplus>0$ and $\mminus>0$.

To state our results, we need some more definitions. 
It is convenient
to define two vectors $\kvec^+=\kvec^+(\musize) = (k_1^+,\ldots, k_m^+)$
and $\kvec^-=\kvec^-(\musize) = (k_1^-,\ldots, k_m^-)$ 
using the following deterministic process:
\begin{itemize}
\item Create a list of length $m$ which contains $\musize(a^+,a^-)$
copies of the ordered pair $(a^+,a^-)$, for all $(a^+,a^-)\in\mathbb{Z}_+^2$;
\item Rearrange this list 
into reverse lexicographic order; 
\item Rename the pairs in this new order as $(k_1^+,k_1^-),\,\, (k_2^+,k_2^-),\ldots, (k_m^+,k_m^-)$.
\end{itemize}
Given any dihypergraph $H=(V,E)\in\calHdra$, there are exactly $\prod_{(a^+,a^-)\in\mathbb{Z}_+^2} \musize(a^+,a^-)!$ ways to order the hyperarcs of $H$ as $e_1,e_2,\ldots, e_m$ such that
$k^+_j=|e_j^+|$ and $k_j^-=|e_j^-|$ for all $j\in [m]$.
By a slight abuse of notation, by fixing one of these orderings we will write $k^+_e$ and $k^-_e$ instead of $k^+_j$ and $k^-_j$, where $e=e_j$. As none of our quantities depend on the particular ordering of the hyperarcs, this does not cause any harm.
Using this notation, we can rewrite condition (\ref{eq:degedgesize}) as
	\begin{equation*}
 \mplus =\sum_{v \in V} d_v^+ = \sum_{e\in E} k_e^+,\qquad
 \mminus=\sum_{v \in V} d_v^- = \sum_{e\in E} k_e^-.
	\end{equation*}
Define
\[ \dmax^+ := \max_{v\in V} d_v^+,\qquad \dmax^- := \max_{v\in V} d_v^-,\qquad
   \kmax^+ := \max_{e\in E} k_e^+,\qquad \kmax^- := \max_{e\in E} k_e^-
   \]
and let
\[ \Delta := \max\{\dmax^+,\dmax^-,\kmax^+,\kmax^-\}.\] 
Write $\dvec=(\dvec^+,\dvec^-)$ and $\kvec = (\kvec^+,\kvec^-)$
and define
\begin{equation}
    \label{eq:R-def}
 R(\dvec) = \sum_{v\in V} d_v^+ d_v^-,\qquad
  R(\kvec) = \sum_{e\in E} k_e^+ k_e^-.
  \end{equation}
Let $\kappa$ be the minimum hyperarc size, defined by
\[ \kappa := \min\big\{ a^+ + a^- \, :\, (a^+,a^-)\in \mathbb{Z}_+^2,\, \musize(a^+,a^-)>0\big\} = \min \big\{ k_e^+ + k_e^-\, :\, e\in E\big\}. \]
We will focus on dihypergraphs with $\kappa\geq 3$, that is, every hyperarc contains at least three vertices. 
Asymptotics are as $n\to\infty$, and we assume that 
$\mplus = \mplus(n)$ and $\mminus = \mminus(n)$ are positive integers which satisfy
\[
\mplus(n)\to\infty \quad \text{ and } \quad \mminus(n)\to\infty \quad \text{ as } \, n\to\infty.
\]

\medskip

Let $(a)_b:= a(a-1)\cdots (a-b+1)$ denote the falling factorial, for all $a\in\mathbb{R}$ and $b\in \mathbb{Z}_+$.
Given two nonnegative integer vectors $\svec=(s_1,\ldots, s_p)$ and $\tvec=(t_1,\ldots, t_q)$ which both sum to $S$, let $S_r:=\sum_{i\in [p]} (s_i)_r$ and $T_r := \sum_{j\in [q]} (t_j)_r$ for $r=2,3$
and define
\begin{equation}
    \label{eq:F-def}
F(\svec,\tvec):= - \frac{S_2T_2}{2S^2} - \frac{S_2T_2}{2S^3} + \frac{S_3T_3}{3S^3} - \frac{S_2T_2(S_2+T_2)}{4S^4} - \frac{S_2^{\, 2}T_3 + S_3 T_2^{\, 2}}{2S^4} + \frac{S_2^{\, 2}T_2^{\, 2}}{2S^5}.
\end{equation} 
This function $F$ arises
from the asymptotic enumeration formula for sparse bipartite graphs with given degrees (see Theorem~\ref{thm:GMW}).

We now state our main result.

\begin{theorem} \label{thm:main}
Suppose that $(\dvec^+,\dvec^-)$ and $\musize$ satisfy \emph{(\ref{eq:degedgesize})}, and that 
$\kappa=\kappa(\musize)\geq 3$. 
Define 
\[ \eta = \eta(\dvec^+,\dvec^-,\musize):= \frac{\Delta^2}{\min\{M^+,M^-\}} + \frac{\Delta^4}{\max\{M^+,M^-\}} + \frac{(\dmax^+ \kmax^+)^3}{(M^+)^2} + \frac{(\dmax^- \kmax^-)^3}{(M^-)^2} \]
If $\eta = o(1)$ then the number of directed hypergraphs with parameters $(\dvec^+,\dvec^-,\musize)$ satisfies
 \begin{align*}	
 \Hrda &= \frac{\mplus!\, \mminus!}{\prod_{v\in V} d_v^+!\,  d_v^-!\, \prod_{e\in E} k_e^+!\, k_e^-!\, \prod_{(a^+,a^-)\in\mathbb{Z}_+^2} \musize(a^+,a^-)!}\\
 & \qquad \times
 \exp\left( F(\dvec^+,\kvec^+) + F(\dvec^-,\kvec^-) - \frac{R(\dvec) R(\kvec)}{M^+\, M^-} + O(\eta) \right).
 \end{align*}
\end{theorem}

\begin{remark} \label{rem:asymptotics}
\emph{In this remark we wish to clarify our use of this asymptotic notation. 
The assumptions of Theorem~\ref{thm:main} can be rewritten in the form $\varepsilon_n= o(1)$ for some expression $\varepsilon_n$. The theorem provides an asymptotic formula with an error term which can be written as $O\big(g_n\big)$, where $g_n$ is a positive function of $n$.  Our proof establishes that there exist positive absolute constants $\varepsilon$ and $C$ such that if $\varepsilon_n < \varepsilon$ then the absolute value of the error term is at most $C g_n$.
We use asymptotic notation in this way throughout the paper, and remark that this usage is consistent with~\cite{GMW}.
} 
\end{remark}

\begin{remark} \label{rem:unnecessary}
\emph{
The expression for $\eta$ always contains an unnecessary term. For example, if $M^-\leq M^+$ then the assumption that $\eta=o(1)$ implies that $\Delta^4=o(M^+)$, and hence
\[ \frac{(\dmax^+ \kmax^+)^3}{(M^+)^2} \leq \frac{\Delta^6}{(M^+)^2} = O\left(\frac{\Delta^4}{M^+}\right).\]
However, it is convenient to write $\eta$ in this symmetric form, to avoid having two different definitions depending on which of $M^+$, $M^-$ is larger.
}
\end{remark}

The following corollary considers the special case of B-hypergraphs and F-hypergraphs.

\begin{corollary}
\label{cor:B-hypergraphs}
Suppose that $(\dvec^+,\dvec^-)$ and $\musize$ satisfy \emph{(\ref{eq:degedgesize})} and 
$\kappa=\kappa(\musize)\geq 3$.    
\begin{itemize}
\item[\emph{(a)}] \emph{[B-hypergraphs]:}\ Suppose that 
$\musize(a^+,a^-)=0$ whenever $a^-\neq 1$. If $\Delta^4=o(M^+)$ then
 \begin{align*}	
&\Hrda \\
 &= 
 \frac{\mplus!\, m!}
 {\prod_{v\in V} d_v^+!\,  d_v^-!\, \prod_{e\in E} k_e^+!\, \prod_{a^+\in\mathbb{Z}_+} \musize(a^+,1)!}\,\exp\left( F(\dvec^+,\kvec^+) 
   - \frac{R(\dvec)}{m} + O\left(\frac{\Delta^4}{M^+}\right)\right).
 \end{align*}
\item[\emph{(b)}] \emph{[F-hypergraphs]:}\  Suppose that
$\musize(a^+,a^-)=0$ whenever $a^+\neq 1$. If $\Delta =o(M^-)$ then
 \begin{align*}	
&\Hrda \\
 &= 
 \frac{\mminus!\, m!}
 {\prod_{v\in V} d_v^+!\,  d_v^-!\, \prod_{e\in E} k_e^-!\, \prod_{a^-\in\mathbb{Z}_+} \musize(1,a^-)!}\, \exp\left( F(\dvec^-,\kvec^-)  
  - \frac{R(\dvec)}{m} + O\left(\frac{\Delta^4}{M^-}\right)\right).
 \end{align*}
\end{itemize}
\end{corollary}

\begin{proof}
For (a), 
the result follows from Theorem~\ref{thm:main} after observing that  $\min\{M^+,M^-\} = m$, $M^+\leq \Delta m$ and $F(\dvec^-,\kvec^-)=0$, and recalling Remark~\ref{rem:unnecessary}. 
The result for (b) follows by symmetry.
\end{proof}

\subsection{Our approach}\label{ss:approach}

In order to make a connection between directed hypergraphs and directed bipartite graphs,
we start by introducing \emph{hyperarc-labelled directed hypergraphs}. 
These are
directed hypergraphs where the hyperarcs are labelled $e_1,e_2,\ldots, e_m$.  
Let $L(\dvec^+,\kvec^+,\dvec^-,\kvec^-)$ denote the number of hyperarc-labelled directed hypergraphs with degree sequence $\dvec:=(\dvec^+,\dvec^-)$ and hyperarc size sequence $\kvec:=(\kvec^+,\kvec^-)$. 
Observe that
\begin{equation}\label{eq:L-H-conversion}
	\frac{L(\dvec^+,\kvec^+,\dvec^-,\kvec^-)}{\Hrda}=\prod_{(a^+,a^-)\in \mathbb{Z}_+^2} \musize(a^+,a^-)!\,,
\end{equation}
since the vectors $\kvec^+$, $\kvec^-$ are fixed, so the only hyperarcs that can be permuted are those
with the same value of $(|e^+|,|e^-|)$.

Next we will exploit the connection between dihypergraphs and directed bipartite graphs. 
With a slight abuse of notation, a hyperarc-labelled directed hypergraph $H=(V,E)$ 
can be represented
by a directed bipartite graph $G=(V,E,W)$, defined as follows. The vertex set of $G$ is $V\cup E$, where one part contains the vertices in $V$ and the other part contains vertices corresponding to the (labelled) elements of $E$. 
The arc set $W$ is defined as follows: if $v\in e^+$ then $W$ contains the arc from $v$ to $e$, while if $v\in e^-$ then $W$ contains the arc from $e$ to $v$.
 This construction is illustrated
in Figure~\ref{f:example} (left to right).

\begin{figure}[ht!]
\begin{center}
\begin{tikzpicture}
\begin{scope}[scale=1.1]
\draw [very thick, ->, rounded corners] (0,3) -- (3.8,3);  
\draw [very thick, ->, rounded corners] (0,3) -- (1.5,2) -- (2.5,2) -- (3.8,2.8); 
\draw [very thick, ->, rounded corners] (1,1) -- (1.5,2) -- (1.6,2)    (2.4,2) -- (2.5,2) -- (3.0,1.2); 
\draw [very thick, ->, rounded corners] (3,1) -- (4,1) -- (4,2.8); 
\draw [very thick, ->, rounded corners] (3.8,1) -- (4,1) -- (4,0.2); 
\draw [very thick, ->, rounded corners] (1,1) -- (0,0) -- (0,2.8); 
\draw [very thick, -, rounded corners] (3,1) -- (0,0) -- (0,0.1); 
\draw [very thick, -, rounded corners] (4,0) -- (0,0) -- (0,0.1); 
\draw [fill] (0,3) circle (0.08);
\draw [fill] (4,3) circle (0.08);
\draw [fill] (1,1) circle (0.08);
\draw [fill] (3,1) circle (0.08);
\draw [fill] (4,0) circle (0.08);
\node [left] at (0,3) {$v_1$};
\node [right] at (4,3) {$v_2$};
\node [left] at (1,1.1) {$v_3$};
\node [right] at (4,0) {$v_4$};
\node [below] at (3.1,0.9) {$v_5$};
\node [above] at (2,3) {$e_4$};
\node [above] at (2,2) {$e_2$};
\node [left] at (0,0) {$e_1$};
\node [right] at (4,1) {$e_3$};
\node [below] at (2,-0.2) {$H$};
\end{scope}
\draw [line width=3pt,<->] (5.75,1.5) -- (7.25,1.5);
\begin{scope}[shift={(8,0)}]
\draw [very thick, ->] (0,3) -- (1.6,3); 
\draw [very thick, ->] (2,3) -- (3.8,3);  
\draw [very thick, ->] (0,3) -- (1.6,2.1); 
\draw [very thick, ->] (2,2) -- (3.8,2.75); 
\draw [very thick, ->] (1,1) -- (1.75,1.75); 
\draw [very thick, ->] (2,2) -- (2.9,1.1); 
\draw [very thick, ->] (3,1) -- (3.6,1); 
\draw [very thick, ->] (4,1) -- (4,2.8); 
\draw [very thick, ->] (4,1) -- (4,0.2); 
\draw [very thick, ->] (1,1) -- (0.25,0.3); 
\draw [very thick, ->] (0,0) -- (0,2.8); 
\draw [very thick, ->] (3,1) -- (0.4,0.2); 
\draw [very thick, ->] (4,0) -- (0.4,0); 
\draw [fill] (0,3) circle (0.08); 
\draw [fill] (4,3) circle (0.08); 
\draw [fill] (1,1) circle (0.08); 
\draw [fill] (3,1) circle (0.08); 
\draw [fill] (4,0) circle (0.08); 
\draw [fill=white] (0,0) circle (0.3); 
\draw [fill=white] (2,2) circle (0.3); 
\draw [fill=white] (4,1) circle (0.3); 
\draw [fill=white] (2,3) circle (0.3); 
\node [left] at (0,3) {$v_1$};
\node [right] at (4,3) {$v_2$};
\node [left] at (1,1.1) {$v_3$};
\node [right] at (4,0) {$v_4$};
\node [below] at (3.1,0.9) {$v_5$};
\node  at (2,3) {$e_4$};
\node  at (2,2) {$e_2$};
\node  at (0,0) {$e_1$};
\node  at (4,1) {$e_3$};
\node [below] at (2,-0.2) {$G$};
\end{scope}
\end{tikzpicture}
\caption{The dihypergraph $H$ (left) corresponds to the directed bipartite graph $G$ (right).}
\label{f:example}
\end{center}
\end{figure}

In this example, $H$ has out-degrees, in-degrees, tail and head sizes given by 
\[ \dvec^+ = (2,0,2,1,2),\quad \dvec^- = (1,3,0,1,1), \quad \kvec^+ = (3,2,1,1),\quad \kvec^- = (1,2,2,1)\]
respectively.  The corresponding directed bipartite graph $G$ has bipartite out-degree sequence
$(\dvec^+,\kvec^-)$ and bipartite in-degree sequence $(\dvec^-,\kvec^+)$.

Equivalently, we can partition $W$ into $W^+$ and $W^-$, where $W^+$ consists of all arcs in $W$ from $V$ to $E$,
and $W^{-}$ consists of all arcs from $E$ to $V$.  After ``forgetting'' the directions of elements of $W^+$ and $W^-$ we see that $G$ corresponds to a pair of (undirected) bipartite graphs $(G^+,G^-)$ where $G^+ = (V,E,W^+)$ has bipartite degree sequence $(\dvec^+,\kvec^+)$
and $G^- = (V,E,W^-)$ has bipartite degree sequence $(\dvec^-,\kvec^-)$.
We will use the terminology ``hyperarc'' for dihypergraphs, ``arc'' for directed bipartite graphs and ``edge'' for bipartite graphs.

For $e\in E$ we write $N_G^-(e)$ for the in-neighbourhood of $e$ in $G$, and $N_G^+(e)$ for the
out-neighbourhood of $e$ in $G$. 
By definition, in any dihypergraph $H$ we have $e^+\cap e^- = \emptyset$ for all $e\in E$, and 
there are no repeated hyperarcs. For the corresponding directed bipartite graph
$G=(V,E,W)$, these conditions state that $G$ has no directed 2-cycles and
\begin{equation}\label{eq:no-repeated-edges}
    \text{ for all distinct $e,f\in E$},\,\, 
 \text{ either } \,\, N^+_G(e) \neq N_G^+(f) \,\, \text{ or } \,\,
    N_G^-(e) \neq N^-_G(f) .
\end{equation}
That is, no two vertices $e,f$
in $E$ have the same in-neighbours and the same out-neighbours in $G$.

Each directed bipartite graph $G$ with no directed 2-cycles
which satisfies (\ref{eq:no-repeated-edges})
uniquely defines a hyperarc-labelled directed hypergraph by letting
\[ e^+ = N^-_{G}(e) \qquad \text{ and } \qquad e^- = N^+_{G}(e) \]
for all $e\in E$.  See Figure~\ref{f:example}, reading from right to left.

Let 
$\calB(\dvec^+,\kvec^+,\dvec^-,\kvec^-)$ be the set of directed bipartite graphs
$G$ such that $G$ has bipartite out-degree sequence $(\dvec^+,\kvec^-)$ and
bipartite in-degree sequence $(\dvec^-,\kvec^+)$. 
Denote by $\calB_0(\dvec^+,\kvec^+,\dvec^-,\kvec^-)$ the set of all
directed graphs in $\calB(\dvec^+,\kvec^+,\dvec^-,\kvec^-)$ with no directed
2-cycles, and write $\B_0(\dvec^+,\kvec^+,\dvec^-,\kvec^-) = |\calB_0(\dvec^+,\kvec^+,\dvec^-,\kvec^-)|$.
Denote by $P(\dvec^+,\kvec^+,\dvec^-,\kvec^-)$ the probability that
a uniformly chosen $G\in \calB_0(\dvec^+,\kvec^+,\dvec^-,\kvec^-)$
satisfies (\ref{eq:no-repeated-edges}).  
Then 
\begin{equation}
L(\dvec^+,\kvec^+,\dvec^-,\kvec^-) = \B_0(\dvec^+,\kvec^+,\dvec^-,\kvec^-)\, 
   P(\dvec^+,\kvec^+,\dvec^-,\kvec^-).
\label{eq:L-B}
\end{equation}

This leads to a strategy for approximating $H(\dvec^+,\dvec^-,\musize)$
which applies whenever it is unlikely that (\ref{eq:no-repeated-edges})
fails (that is, when $P(\dvec^+,\kvec^+,\dvec^-,\kvec^-)$ is close to 1):
\begin{itemize}
\item[(i)] Define vectors $\kvec^+,\kvec^-$ as described in Section~\ref{s:statements}.
\item[(ii)] 
Use a switching argument to obtain an asymptotic formula for $\B_0(\dvec^+,\kvec^+,\dvec^-,\kvec^-)$.
\item[(iii)] Prove that only a negligible proportion of these directed bipartite graphs $G$
fail (\ref{eq:no-repeated-edges}). 
\item[(iv)]  Apply (\ref{eq:L-H-conversion}) and (\ref{eq:L-B}) to deduce an 
asymptotic enumeration formula for $H(\dvec^+,\dvec^-,\musize)$.
\end{itemize}
We remark that in order to have a repeated hyperarc, there must be at least
one pair $(a^+,a^-)$ with $\musize(a^+,a^-)\geq 2$.  If $\musize$ only
takes values in $\{0,1\}$ then $P(\dvec^+,\kvec^+,\dvec^-,\kvec^-)=1$
and step~(iii) is unnecessary.

\bigskip

To conclude this section we outline the structure of the rest
of this paper.
In Section~\ref{s:bipartite} we introduce our notation and some asymptotic enumeration
results for sparse bipartite graphs which will be useful in our proof.
In Section~\ref{s:switching} we apply the switching method to obtain
an asymptotic enumeration
formula for $\B_0(\dvec^+,\kvec^+,\dvec^-,\kvec^-)$
under the assumptions of Theorem~\ref{thm:main}.
Finally in Section~\ref{s:repeats} we complete the proof of Theorem~\ref{thm:main} by showing that
$P(\dvec^+,\kvec^+,\dvec^-,\kvec^-)$ is close to~1 under the
assumptions of the theorem.

\section{Preliminaries}\label{s:bipartite}

Let $V$ and $W$ be disjoint finite sets.
We work with bipartite degree sequences $(\svec,\tvec)$ of nonnegative integers such that $s_v$ is the degree of vertex $v$ for all $v\in V$ and $t_w$ is the degree of vertex $W$ for all $w\in W$, and $S = \sum_{v\in V} s_v = \sum_{w\in W} t_w$.
Denote by  $\smax,\tmax$ the maximum component of the vectors $\svec,\tvec$ respectively, and define $F(\svec,\tvec)$ as in (\ref{eq:F-def}).

Let $\mathcal{B}(\svec,\tvec)$ be the set of bipartite graphs with vertex bipartition $V\cup W$ 
and bipartite degree sequence $(\svec,\tvec)$. 
We will use the asymptotic enumeration result for $|\mathcal{B}(\svec,\tvec)|$,
restated from~\cite{GMW}. 

\begin{theorem}[{\cite[Theorem~1.3]{GMW}}] 
\label{thm:GMW}
Suppose that $S\to\infty$, and that $1\leq \smax\tmax = o(S^{2/3})$. Then
\begin{align*} 
|\mathcal{B}(\svec,\tvec)| = \frac{S!}{\prod_{v\in V} s_v! \prod_{w\in W} t_w!}
 \, \exp\left( F(\svec,\tvec) + O\left(\frac{\smax^3\tmax^3}{S^2}\right)\right).
 \end{align*}
\end{theorem}

Note that there are other asymptotic enumeration results for bipartite graphs by degree sequence which hold under less restrictive conditions on the degree sequence, including~\cite{GM09} in the dense regime, and most recently by Liebenau and Wormald~\cite{LW} for a wide range of near-regular degrees. However, the switching argument in Section~\ref{s:switching} imposes an upper bound on the maximum degrees and hyperarc sizes, which implies that the result of Theorem~\ref{thm:GMW} is suitable for our purposes. 

We can apply this result to  write down an asymptotic enumeration formula for
the number of sparse directed bipartite graphs with given degrees, using the fact that
\[ \B(\dvec^+,\kvec^+,\dvec^-,\kvec^-) = |\mathcal{B}(\dvec^+,\kvec^+)|\, |\mathcal{B}(\dvec^-,\kvec^-)|.\]

\begin{corollary}
\label{cor:size-B}
Suppose that $(\dmax^+ \, \kmax^+)^{3/2} = o(M^+)$ and $(\dmax^-\,\kmax^-)^{3/2} = o(M^-)$.  
Then
\begin{align*}
 \B(\dvec^+,\kvec^+,& \dvec^-,\kvec^-)\\
 &=
 \frac{\mplus!\, \mminus!}
 {\prod_{v\in V} d_v^+!\,  d_v^-!\, \prod_{e\in E} k_e^+!\, k_e^-!}
 \\ & \quad \times 
  \exp\left( F(\dvec^+,\kvec^+) + F(\dvec^-,\kvec^-)  + 
  O\left( \frac{(\dmax^+ \kmax^+)^3}{(M^+)^2} + \frac{(\dmax^- \kmax^-)^3}{(M^-)^2} 
  \right)\right).
 \end{align*}
 \end{corollary}

We will use the following result of McKay~\cite{McKay1981} to obtain upper bounds on the probability of various bad cases occuring. In fact we use a restatement of a special case of McKay's result, using the wording from~\cite[Lemma~2.1]{BG2016b}.

\begin{theorem}[{\cite[Theorem 3.5(a)]{McKay1981}}]
Let $\mathcal{B}(g)$ denote the set of bipartite graphs with vertex bipartition given by $\{a_1,\ldots,a_n\}\cup\{b_1,\ldots,b_m\}$ and degree sequence
\[
\gvec = (g_1,\ldots,g_n; g_1',\ldots,g_m').
\]
(Here vertex $a_i$ has degree $g_i$ for $i=1,\ldots,n$, and vertex $b_j$ has degree $g_j'$ for $j=1,\ldots,m$.)  
Let $L$ be a subgraph of the complete bipartite graph on this vertex bipartition, and let $\mathcal{B}(\gvec,L)$ be the set of bipartite graphs in $\mathcal{B}(\gvec)$ which contain $L$ as a subgraph.  
Write $E_{\gvec} = \sum_{i=1}^n g_i$ and $E_{\ellvec} = \sum_{i=1}^n \ell_i$, where $\ellvec = (\ell_1,\ldots,\ell_n;\ell_1',\ldots,\ell_m')$ is the degree sequence of $L$.  
Finally, let $g_{\max}$ and $\ell_{\max}$ denote the maximum degree in $g$ and $\ell$, respectively, and define
\[
\Gamma = 2g_{\max}(g_{\max} + \ell_{\max} - 1) + 2.
\]
If $E_{\gvec} - \Gamma \geq E_{\ellvec}$ then
\[
\frac{|\mathcal{B}(\gvec,L)|}{|\mathcal{B}(\gvec)|} \leq 
\frac{\prod_{i=1}^n (g_i)_{\ell_i} \prod_{j=1}^m (g_j')_{\ell_j'}}{(E_{\gvec} - \Gamma)_{E_{\ellvec}}}.
\]
\label{thm:subgraph}
\end{theorem}

\section{Probability of no directed 2-cycles}\label{s:switching}

Recall the definition of $R(\dvec)$, $R(\kvec)$ from (\ref{eq:R-def}).
Observe that 
\begin{equation}
\label{eq:Rbounds}
\left.
\begin{aligned}
 R(\dvec) &\leq \min\{ \dmax^+ M^-,\, \dmax^- M^+\} \leq \Delta\min\{ M^+,M^-\},\\
\max\{M^+,M^-\}\leq R(\kvec)&\leq \min\{ \kmax^+ M^-,\, \kmax^- M^+\} \leq \Delta\min\{M^+,M^-\},
\end{aligned}
\qquad \right\}
\end{equation}
using the assumption that $k_e^+,k_e^-\geq 1$ for all $e\in E$ for the lower bound on $R(\kvec)$.
Define 
\[ N = \max\{\lceil \log\big(\min\{M^+,\, M^-\}\big)\rceil,\, \lceil 20\mu\rceil\}\] 
where 
\[ \mu:= \frac{R(\dvec)\, R(\kvec)}{M^+ M^-}.\]
For $t=0,\ldots, N$ let $\mathcal{S}_t$ be the set of all directed bipartite graphs in
$\calB(\dvec^+,\kvec^+,\dvec^-,\kvec^-)$ which contain exactly $t$ directed 2-cycles.
First we prove that with high probability there are at most $N$ directed 2-cycles in a 
randomly chosen element of $\calB(\dvec^+,\kvec^+,\dvec^-,\kvec^-)$.

\begin{lemma}
\label{lem:cutoff-N}
Suppose that $\Delta^2  = o\big(\min\{ M^+,M^-\}$ and $\Delta^4 = o(\max\{M^+,M^-\})$.  Then
\[ \sum_{t=0}^N |\mathcal{S}_t| = \left(1 + O\left(\frac{\Delta^2}{\min\{M^+,\, M^-\}}\right)\right)\, \vec{B}(\dvec^+,\kvec^+,\dvec^-,\kvec^-).\]
\end{lemma}

\begin{proof}
For ease of notation, let $Q:= N+1$, and let $Z$ be the number of sets of $Q$ distinct directed 2-cycles in a uniformly chosen
element of $\calB(\dvec^+,\kvec^+,\dvec^-,\kvec^-)$.  Equivalently, we can 
let $Z$ be the number of sets of $Q$ distinct pairs $(v,e)\in V\times E$ which are edges in both $G^+$ and $G^-$,
where $G^+$ is a uniformly random bipartite graph in $\mathcal{B}(\dvec^+,\kvec^+)$
and $G^-$ is a uniformly random bipartite graph in $\mathcal{B}(\dvec^-,\kvec^-)$. 
We find an upper bound on $\E Z$ using Theorem~\ref{thm:subgraph}.

Let $A$ be a set of $Q$ distinct pairs $(v_1,e_1),\ldots, (v_Q,e_Q)\in V\times E$. 
Applying Theorem~\ref{thm:subgraph} to both $G^+$ and $G^-$, the probability that both $G^+$ and $G^-$ contain all edges of $A$ is at most
\begin{align}
 \frac{\prod_{ve\in A} d_v^+ d_v^- k_e^+ k_e^-}{\big(M^+ - O\big(\Delta^2\big)\big)_{Q}\big(M^- - O\big(\Delta^2\big)\big)_{Q}}
 \leq   \frac{ \prod_{ve\in A} d_v^+ d_v^- k_e^+ k_e^-}{\big(M^+M^-\big)^{Q}}\, \left(1 + O\left(\frac{Q\, \Delta^2}{\min\{ M^+,M^-\}}\right)\right), \label{eq:attempt}
\end{align}
as $Q=O(\Delta^2)$.
Now we must sum over all choices of $A$. 
We will choose the edges one by one in a sequence, and divide by $Q!$ later to produce a set. Supposing that $(v_1,e_1),\ldots, (v_{i-1},e_{i-1})$ have been chosen,
the contribution to the numerator of (\ref{eq:attempt}) from $(v_i,e_i)$ is at most
\[  \sum_{v\in V}
 d_v^+ d_v^- \sum_{e\in E} k_e^+ k_e^-
 =  R(\dvec) R(\kvec).
\]
Combining this with (\ref{eq:attempt}) we find that \
 \begin{align*}
\E Z \leq \frac{1}{Q!}\, \left(\frac{ R(\dvec)R(\kvec) }{M^+M^-}\right)^{Q} \left(1 + O\left(\frac{Q\, \Delta^2}{\min\{ M^+,M^-\}}\right)\right)
  = \frac{\mu^Q}{Q!}  \left(1 + O\left(\frac{Q\, \Delta^2}{\min\{ M^+,M^-\}} \right)\right).
\end{align*}
  We divide by $Q!$ as $Z$ counts \emph{sets} of $Q$ distinct directed 2-cycles, while our calculations above involved \emph{sequences}. 
Now
\[ \frac{\mu^Q}{Q!}\leq \left(\frac{ e \mu}{Q}\right)^{Q} \leq \left(\frac{ e}{8}\right)^{Q} \leq \frac{1}{\min\{M^+,\, M^-\}},\]
using the fact that $Q\geq 20\mu$ and $Q\geq \log(\min\{M^+,M^-\})$.
Furthermore, our assumptions imply that $\Delta=o(\sqrt{\min\{M^+,M^-\}})$ and hence
\begin{align*} 
\frac{Q \, \Delta^2}{\min\{ M^+,M^-\}}
  &\leq \frac{\Delta^2\big(\Delta^2 +  \log(\min\{M^+,M^-\})\big)}{\min\{ M^+,M^-\}}
  = O(\Delta^2).
 \end{align*}
Therefore by Markov's lemma, the probability that $\vec{B}(\dvec^+,\kvec^+,\dvec^-,\kvec^-)$ contains more than $N$ distinct directed 2-cycles is at most $\E Z = O\big(\Delta^2/\min\{M^+,M^-\}\big)$. The result follows.
\end{proof}

\bigskip

Note that we are interested in the set $\calB_0(\dvec^+,\kvec^+,\dvec^-,\kvec^-) = \mathcal{S}_0$.  We will obtain an asymptotic formula for $|\mathcal{S}_0|$ using the
switching method. The following \emph{forward switching} is designed to reduce the number of directed 2-cycles by exactly one, and is illustrated in Figure~\ref{fig:switching}. 
The forward switching operation proceeds as follows.  
Starting from some directed bipartite graph
$G\in \mathcal{S}_t$:
\begin{itemize}
\item Choose a directed 2-cycle on vertices $v,e$ and two other arcs $w_1f_1$ and $f_2w_2$, where $v,w_1,w_2\in V$  and $e, f_1,f_2\in E$, such that $w_1\neq w_2$, $f_1\neq f_2$ and 
\[ \{ f_1 w_1, w_2f_2,\, \, f_1v,\, vf_1,\, w_1 e,\, ew_1,\, w_2 e,\, ew_2,\, f_2 v,\, vf_2\}\cap E(G) = \emptyset.\]
\item  Define the directed bipartite graph $G'$ from $G$ by deleting the arcs $w_1f_1, ve, ev, f_2w_2$ and inserting the arcs $w_1e, ew_2, f_2v, vf_1$.
\end{itemize}
The conditions in the first step ensure that $G'\in\mathcal{S}_{t-1}$.  Furthermore, in a switching we cannot have $w_1=v$ as $ve\in E(G)$ and $w_1e\not\in E(G)$.  Similarly
the conditions of the switching imply that $w_2\neq v$, $f_1\neq e$ and $f_2\neq e$.

\begin{figure}[ht!]
\begin{center}
\begin{tikzpicture}
\draw [->,very thick] (0,2) to[out=-80,in=80] (0,0.2);  
\draw [->,dashed,black!40] (-0.12,0.15) to[out=100,in=-100] (-0.1,1.85); 
\draw [->,very thick] (4,0) to[out=100,in=-100] (4,1.8);  
\draw [->,dashed,black!40] (4.12,2) to[out=-80,in=80] (4.1,0.15);  
\draw [->, very thick] (2,0) to[out=80,in=-80] (2.1,1.8);  
\draw [->, very thick] (2,2) to[out=-100,in=100] (1.9,0.2);  
\draw [->,dashed,black!40] (0,2) -- (1.85,-0.1);  
\draw [->,dashed,black!40] (1.85,0.15) to[out=125,in=-35] (0.15,2.1);  
\draw [->,dashed,black!40] (2,2) -- (0.15,-0.1);  
\draw [->,dashed,black!40] (0.15,0.15) to[out=55,in=-145] (1.85,2.1);  
\draw [->,dashed,black!40] (2,-0.1) -- (3.85,2.1);  
\draw [->,dashed,black!40] (3.85,1.85) to[out=-130,in=35] (2.15,-0.1);  
\draw [->,dashed,black!40] (4,-0.1) -- (2.15,2.1);  
\draw [->,dashed,black!40] (2.15,1.9) to[out=-55,in=145] (3.85,-0.1);  
\draw [fill] (0,0) circle (0.1);  \draw [fill] (2,0) circle (0.1); \draw [fill] (4,0) circle (0.1);
\draw [fill] (0,2) circle (0.1);  \draw [fill] (2,2) circle (0.1); \draw [fill] (4,2) circle (0.1);
\node [below] at (0,-0.1) {$f_1$}; \node [below] at (2,-0.1) {$e$}; \node [below] at (4,-0.1) {$f_2$};
\node [above] at (0,2.1) {$w_1$}; \node [above] at (2,2.1) {$v$}; \node [above] at (4,2.1) {$w_2$};
\node [below] at (2,-0.8) {$G\in \mathcal{S}_{t}$};
\draw [->,line width=4pt] (5,1) -- (6,1);
\begin{scope}[shift={(7,0)}]
\draw [->,dashed,black!40] (0,2) to[out=-80,in=80] (0,0.2);  
\draw [->,dashed,black!40] (-0.12,0.15) to[out=100,in=-100] (-0.1,1.85); 
\draw [->,dashed,black!40] (4,0) to[out=100,in=-100] (4,1.8);  
\draw [->,dashed,black!40] (4.12,2) to[out=-80,in=80] (4.1,0.15);  
\draw [->, dashed,black!40] (2,0) to[out=80,in=-80] (2.1,1.8);  
\draw [->, dashed,black!40] (2,2) to[out=-100,in=100] (1.9,0.2);  
\draw [->,very thick] (0,2) -- (1.85,-0.1);  
\draw [->,dashed,black!40] (1.85,0.15) to[out=125,in=-35] (0.15,2.1);  
\draw [->,very thick] (2,2) -- (0.15,-0.1);  
\draw [->,dashed,black!40] (0.15,0.15) to[out=55,in=-145] (1.85,2.1);  
\draw [->,very thick] (2,0) -- (3.85,2.1);  
\draw [->,dashed,black!40] (3.85,1.85) to[out=-130,in=35] (2.15,-0.1);  
\draw [->,very thick] (4,0) -- (2.15,2.1);  
\draw [->,dashed,black!40] (2.15,1.9) to[out=-55,in=145] (3.85,-0.1);  
\draw [fill] (0,0) circle (0.1);  \draw [fill] (2,0) circle (0.1); \draw [fill] (4,0) circle (0.1);
\draw [fill] (0,2) circle (0.1);  \draw [fill] (2,2) circle (0.1); \draw [fill] (4,2) circle (0.1);
\node [below] at (0,-0.1) {$f_1$}; \node [below] at (2,-0.1) {$e$}; \node [below] at (4,-0.1) {$f_2$};
\node [above] at (0,2.1) {$w_1$}; \node [above] at (2,2.1) {$v$}; \node [above] at (4,2.1) {$w_2$};
\node [below] at (2,-0.8) {$G'\in \mathcal{S}_{t-1}$};
\end{scope}
\end{tikzpicture}
\caption{The switching (from left to right) and reverse switching (from right to left)}
\label{fig:switching}
\end{center}
\end{figure}

A \emph{reverse switching} is the reverse of the forward switching, and proceeds as follows. Starting with a bipartite directed graph $G'\in\mathcal{S}_{t-1}$,
\begin{itemize}
\item Choose a directed 2-path $w_1 e w_2$ starting and ending in $V$, and a directed 2-path $f_2 v f_1$ starting and ending in $E$, such that 
\[ \{ f_1 w_1, w_1 f_1, e w_1, f_1v, ev, ve, vf_2, w_2 e, w_2f_2, f_2w_2\} \cap E(G') = \emptyset.\]
\item Define the directed bipartite graph $G$ from $G'$ by deleting the arcs $w_1 e$, $e w_2$, $f_2 v$, $v f_1$ and inserting the arcs $w_1f_1, ve, ev, f_2w_2$.
\end{itemize}
The reverse switching is illustrated in Figure~\ref{fig:switching}, reading right-to-left.
The conditions in the first step guarantee that $G\in \mathcal{S}_t$ and that $v\neq w_j$ and $e\neq f_j$ for $j=1,2$.

\begin{lemma}
Assume that $\Delta^2=o(\min\{M^+,M^-\})$. 
Then 
\[ |\mathcal{S}_t| = \frac{R(\dvec)R(\kvec) + O\left( \Delta^3 \min\{R(\dvec),R(\kvec)\} + (t-1)\Delta \big(R(\dvec) + R(\kvec)\big)\right)}{t(M^+-t)(M^- -t)} \, |\mathcal{S}_{t-1}|\]
uniformly for all $t\in [N]$ such that $\mathcal{S}_{t-1}$ is nonempty.
\label{lem:switchings-analysis}
\end{lemma}

\begin{proof}
Let $G\in\mathcal{S}_t$. There are at most $t(M^+ -t) (M^--t)$ ways to perform a forward switching starting from $G$,
as there are exactly $t$ choices for the directed 2-cycle, at most $M^+-t$ ways to
choose the arc $w_1f_1$ so that it is not part of a directed 2-cycle, and at most $M^- - t$ ways to choose the arc $f_2w_2$ so that it is not part of a directed 2-cycle.
For a lower bound, we must subtract an upper bound on the number of choices which do not correspond to a forward switching, because one or more of the conditions are not satisfied.
\begin{itemize}
\item  There are at most $t\bigl(R(\dvec)+ R(\kvec)\bigr) = O\bigl(t\Delta \min\{M^+,M^-\}\bigr)$ choices where $w_1=w_2$ or $f_1=f_2$, using (\ref{eq:Rbounds}).
\item 
There are at most 
\begin{align*}
  &  tM^-\big( 2\dmax^+ \kmax^+ + \dmax^+ \kmax^- + \dmax^-\kmax^+\big) + tM^+ \big( 2\dmax^-\kmax^- + \dmax^+ \kmax^- + \dmax^- \kmax^+\big)\\
    &= O\bigl( t \Delta^2\, (M^+ + M^-)\bigr)
\end{align*}
choices where one of the eight arcs $w_1e$, $ew_1$, $w_2 e$, $e w_2$, $vf_1$, $f_1v$, $vf_2$, $f_2v$ is already present in $G$.
\end{itemize}
After subtracting the number of bad cases and comparing with the upper bound, it follows
that the number of forward switchings from $G$ is
\begin{align}
\label{eq:forward}
 &t \left( (M^+ -t)(M^- - t) - O\bigl( \Delta^2(M^+ + M^-)\bigr)\right) \notag \\
 & {} \qquad = t(M^+-t)(M^- -t)\left( 1 + O\left(\frac{\Delta^2}{\min\{ M^+,M^-\}}\right)\right)
\end{align} 
for all $G\in\mathcal{S}_t$.

Next we fix $G'\in\mathcal{S}_{t-1}$ and analyse the number of \emph{reverse switchings} that can be performed from $G'$ to produce an
element of $\mathcal{S}_t$.  There are at most $R(\dvec)R(\kvec)$ ways to choose the two directed 2-paths.
From this expression we must subtract the number of choices which violate a condition of the reverse switching.
\begin{itemize}
\item  There are at most $(t-1)\bigl(R(\dvec) + R(\kvec)\bigr)$ choices with $w_1=w_2$ or $f_1=f_2$.
\item The number of choices where all six vertices are distinct but one of the chosen arcs belongs to a directed 2-cycle in $G'$ (that is, one of the arcs $ew_1,f_1v,vf_2,w_2e$ is present) is at most
\[ (t-1)\left( (\kmax^- + \kmax^+) R(\dvec)
  + (\dmax^- + \dmax^+) R(\kvec)\right).\]
\item  
Assume first that $R(\dvec) \leq R(\kvec)$.  
Then there are at most
\begin{align*}
 \left( 2(\dmax^+ + \dmax^-)\kmax^+\kmax^- + \dmax^-(\kmax^+)^2 + \dmax^+(\kmax^-)^2\right)\, R(\dvec) = O\bigl( \Delta^3 R(\dvec)\bigr)
\end{align*}
choices where all six vertices are distinct but one of the remaining six arcs $ve$, $ev$, $w_1f_1$, $f_1w_1$, $w_2f_2$, $f_2w_2$ is already present in $G'$. 
A similar bound holds if $R(\kvec) < R(\dvec)$.
Overall we can certainly say that there are at most
\[ O\big( \Delta^3 \min\{ R(\dvec), R(\kvec)\}\big)\]
such choices. 
\end{itemize}
By subtracting these bad choices and comparing with the upper bound, we find that the number of reverse switchings from $G'$ is
\begin{equation}
\label{eq:reverse}
 R(\dvec) R(\kvec) - O\left(  \Delta^3\min\{ R(\dvec), R(\kvec)\} + (t-1)\Delta \bigl(R(\dvec)+R(\kvec)\bigr) \right).
 \end{equation} 
The lemma follows from (\ref{eq:forward}) and (\ref{eq:reverse}) by double-counting the number of switchings from $\mathcal{S}_t$ to $\mathcal{S}_{t-1}$, and applying (\ref{eq:Rbounds}).
\end{proof}

We will combine the estimates from Lemma~\ref{lem:switchings-analysis} using the following summation
lemma.

\begin{lemma}[{\cite[Lemma 2.5]{BG2016a}},{\cite[Corollary 4.5]{GMW}}]
\label{lemma:sum}
Let $N \geq 2$ be an integer and, for $1 \leq i \leq N$, let real numbers $A(i), C(i)$ be given such that $A(i) \geq 0$ and $A(i) - (i-1)C(i) \geq 0$. Define
\[
A_1 = \min_{1 \leq i \leq N} A(i), \quad 
A_2 = \max_{1 \leq i \leq N} A(i), \quad 
C_1 = \min_{1 \leq i \leq N} C(i), \quad 
C_2 = \max_{1 \leq i \leq N} C(i).
\]
Suppose that there exists a real number $\hat{c}$ with $0 < \hat{c} < \tfrac{1}{3}$ such that 
$\max \Bigl\{A_2/N,\, |C_1|,\, |C_2|\Bigr\} \leq \hat{c}$.
Define $n_0, \ldots, n_N$ by $n_0 = 1$ and
\[
n_i = \frac{1}{i}\Bigl(A(i) - (i-1)C(i)\Bigr)n_{i-1}
\]
for $i=1,\ldots, N$. Then
\[
\Sigma_1 \;\leq\; \sum_{i=0}^N n_i \;\leq\; \Sigma_2,
\]
where
\[
\Sigma_1 = \exp\Bigl(A_1 - \tfrac{1}{2}A_1C_2\Bigr) - (2e\hat{c})^N,
\]
\[
\Sigma_2 = \exp\Bigl(A_2 - \tfrac{1}{2}A_2C_1 + \tfrac{1}{2}A_2C_1^2\Bigr) + (2e\hat{c})^N.
\]
\end{lemma}

\bigskip

\begin{lemma}\label{lem:nocycles}
Suppose that 
$\Delta^2 = o\bigl(\min\{M^+,M^-\}\bigr)$  
and 
$\Delta^4 = o\bigl(\max\{M^+,M^-\}\bigr)$.   
Then
\[ \sum_{t=0}^N |\mathcal{S}_t|
= |\mathcal{S}_0|\, \exp\left( \frac{R(\dvec)R(\kvec)}{M^+ M^-} + O\left(
\frac{\Delta^4}{\max\{M^+, M^- \}}\right) \right)\, . \]
\end{lemma}

\begin{proof}
By (\ref{eq:forward}), any $G\in\mathcal{S}_t$ can be transformed into some $G'\in\mathcal{S}_{t-1}$ by a forward switching.  
Therefore, if $\mathcal{S}_0$ is empty then $\mathcal{S}_t$ is empty for all $t=1,\ldots, N$.  The lemma holds in this case, and so we assume that $\mathcal{S}_0\neq \emptyset$.

Define
\[ A_0 = \frac{R(\dvec)R(\kvec)}{(M^+-1)(M^- - 1)}.
\]
Then expanding the left hand side of the expression below around $t=1$ gives
\begin{align*}
    \frac{R(\dvec)R(\kvec)}{(M^+-t)(M^- -t)} 
    &= A_0 + O\left(\frac{(t-1)\, R(\dvec)R(\kvec)(M^+ + M^-)}{(M^+\, M^-)^2}\right)\\
    &= A_0  + O\left(\frac{(t-1)\Delta \big(R(\dvec)+R(\kvec)\big)}{M^+M^-}\right)
\end{align*}
using (\ref{eq:Rbounds}).
Combining this with Lemma~\ref{lem:switchings-analysis}, which applies by our assumptions, we obtain
\[
  |\mathcal{S}_t| = 
  \frac{|\mathcal{S}_{t-1}|}{t}\, \left( A_0  + O\left( \frac{\Delta^3\, \min\{R(\dvec),R(\kvec)\}}{M^+ M^-} + (t-1)\frac{\Delta \big(R(\dvec) + R(\kvec)\big)}{M^+ M^-}
  \right)\right)
  \]
uniformly for any $t\in [N]$ such that $\mathcal{S}_{t-1}$ is nonempty. 
Hence we
can define a real number $\alpha_t$ for all $t=1,\ldots, N$ such that
\begin{equation}
    \label{eq:AC-alpha} 
    |\mathcal{S}_t| = 
  \frac{|\mathcal{S}_{t-1}|}{t}\, \left( A_0 + \frac{\alpha_t \Delta^3\min\{R(\dvec),R(\kvec)\}}{M^+M^-}  + (t-1) \frac{\alpha_t \Delta \big(R(\dvec) + R(\kvec)\big)}{M^+M^-}\right),
  \end{equation}
where $|\alpha_t|$ is bounded independently of $t$ and $n$. In particular, if $\mathcal{S}_{t-1}$ is nonempty then $\alpha_t$ is uniquely defined by (\ref{eq:AC-alpha}), and if $\mathcal{S}_{t-1}$ is empty then we let $\alpha_t=0$.
  Now define
  \begin{align*}
A(t) &:= A_0 + \frac{\alpha_t \Delta^3\min\{R(\dvec),R(\kvec)\}}{M^{+} M^{-} } 
 = \mu+O\left(\frac{\Delta^4}{\max\{M^+,M^-\}}\right)=\mu+o(1),\\
C(t) &:=  - \frac{\alpha_t \Delta \big(R(\dvec) +R(\kvec)\big)}{M^+M^-} = O\left(\frac{\Delta^2}{\max\{M^+,M^-\}}\right) = o(1)
\end{align*}
for $t=1, \ldots, N$.
Then we can rewrite (\ref{eq:AC-alpha}) as
\[ |\mathcal{S}_t| = \frac{|\mathcal{S}_{t-1}|}{t}\, \big( A(t) - (t-1)C(t)\big).\]

We now check the assumptions of Lemma~\ref{lemma:sum}. First we claim that
$A(t) - (t-1)C(t)\geq 0$ for all $t=1,\ldots, N$.  If $\mathcal{S}_{t-1}$ is nonempty then (\ref{eq:AC-alpha}) implies that $A(t)-(t-1)C(t)\geq 0$, since $|\mathcal{S}_t|\geq 0$.
Otherwise $\mathcal{S}_{t-1}$ is empty, and hence $A(t) = A_0 > 0$ and $C(t)=0$,
and so $A(t)-(t-1)C(t) >0$.  This establishes the first claim. 
Next, if $\alpha_t\leq 0$ then $A(t)\geq A(t)-(t-1)C(t)$ and we already know that $A(t)-(t-1)C(t)\geq 0$.  Hence $A(t)\geq 0$ in this case, while if $\alpha_t> 0$ then $A(t)\geq 0$ by definition. 

Next, define $A_1 = \min_{1\leq t\leq N} A(t)$, $A_2 = \max_{1\leq t\leq N} A(t)$, $C_1 = \min_{1\leq t\leq N} C(t)$
and $C_2 = \max_{1\leq t\leq N} C(t)$.  Let $\hat{c} = 1/16$. Since $A(t) = \mu+o(1)$ it follows that $\max\{ A_2/N,\, |C_1|,|C_2|\} \leq \hat{c}$ for $n$ sufficiently large, by definition of $N$.

Hence Lemma~\ref{lemma:sum} applies. 
For $t=1,\ldots, N$, since $A(t)=\mu+O(1)$ we obtain
\begin{align*}
A(t) C_1,\, A(t)C_2 &= O\left(\frac{\Delta^2}{\max\{M^+,M^-\}}\left(\frac{R(\dvec) R(\kvec)\big)}{M^+M^-}+1\right)\right) = O\left(\frac{\Delta^4}{\max\{ M^+,\, M^-\}}\right)
\end{align*} 
using (\ref{eq:Rbounds}).  
Therefore
\[ A_2 - \dfrac{1}{2}A_2 C_1 +\frac{1}{2}A_2C_1^2= 
\mu + O\left( \frac{\Delta^4}{\max\{M^+,M^-\}} 
\right),
\]
and the same expression holds for $A_1 - \dfrac{1}{2} A_1C_2$ within the stated error bounds.
Combining the upper and lower bounds from Lemma~\ref{lemma:sum} gives
\[ \sum_{t=0}^N \frac{|\mathcal{S}_t|}{|\mathcal{S}_0|}  =
  \exp\left(\frac{R(\dvec)R(\kvec)}{M^{+} M^{-} } +  O\left(\frac{\Delta^4}{\max\{M^+, M^- \}}\right)\right) + O\bigl( (e/8)^N\bigr).  \]
Finally, note that the sum we are estimating is at least one, and
\[ (e/8)^N  \leq (e/8)^{\log(\max\{M^+,M^-\})}\leq \frac{1}{\max\{M^+\, M^-\}} < \frac{\Delta^4}{\max\{ M^+,M^-\}}.\]
Therefore we can absorb the additive term into the relative error bound, completing the proof.
\end{proof}

The following Corollary summarizes the results of this section. It is proved by combining Lemma~\ref{lem:cutoff-N} with Lemma~\ref{lem:nocycles}, recalling that $\mathcal{S}_0 = \calB_0(\dvec^+,\kvec^+,\dvec^-,\kvec^-)$. 

\begin{corollary}~\label{cor:B0}
Suppose that 
$\Delta^2 = o\bigl(\min\{M^+,M^-\}\bigr)$ and
$\Delta^4 = o\bigl(\max\{M^+,M^-\}\bigr)$. Then
\begin{align*}
&\B_0(\dvec^+,\kvec^+,\dvec^-,\kvec^-)\\
&=\exp\left(-\frac{R(\dvec) R(\kvec)}{M^+M^-}+O\left(\frac{\Delta^4}{\max\{M^+,M^-\}}+\frac{\Delta^2}{\min\{M^+,M^-\}}\right)\right)\B(\dvec^+,\kvec^+,\dvec^-,\kvec^-).
\end{align*}
\end{corollary}

\section{Repeated hyperarcs}\label{s:repeats}

The final step is to show that under the conditions of our theorems,
there is a very low probability that a directed bipartite graph gives rise
to a hyperarc-labelled dihypergraph with repeated hyperarcs.
Denote by
$Q(\dvec^+,\kvec^+,\dvec^-,\kvec^-)$ the number of directed bipartite graphs
$G\in \calB_0(\dvec^+,\kvec^+,\dvec^-,\kvec^-)$ 
such that $N^+_{G}(e)=N^+_G(f)$ and $N^-_G(e)=N^-_G(f)$; that is, $e$ and $f$ describe the same hyperarc in the corresponding dihypergraph.

If we remove two copies of a repeated hyperarc from a dihypergraph $H$,
then we obtain a dihypergraph $H'$ where the degree sequences of the corresponding
directed bipartite graph $G$ change in the following manner: 
\begin{itemize}
    \item The in-degree and out-degree
of two vertices $e,f$ (corresponding to the two removed hyperarcs of $H$) are now zero. (We could remove these two vertices from $E$ but in our enumeration it is convenient to keep them as isolated vertices.)
\item The out-degree of each vertex in $N_G^+(e)$ is decreased by~2, and the in-degree of each vertex in $N_G^-(e)$ is decreased by~2.
\end{itemize}
We analyse this situation in the following lemma.

\begin{lemma}
\label{lem:nomults}
Let $\dvec^+,\kvec^+,\dvec^-,\kvec^-$ be vectors of positive integers with minimum hyperarc size $\kappa\geq 3$. Then for
$\Delta^2 = o\bigl(\min\{M^+,M^-\}\bigr)$ and 
$\Delta^4 = o\bigl(\max\{M^+,M^-\}\bigr)$ we have
\[
Q(\dvec^+,\kvec^+,\dvec^-,\kvec^-)=\B_0(\dvec^+,\kvec^+,\dvec^-,\kvec^-)\, O\left(\frac{\Delta^3}{\max\{M^+,M^-\}}\right).
\]
\end{lemma}

\begin{proof}
Let $e,f\in E$ with $k_{e}^+=k_{f}^+$ and $k_{e}^-=k_f^{-}$. In addition, let $U_e^+, U_e^- \subseteq V$ be disjoint sets of vertices such that 
\[ |U_e^{+}|=k_e^+, \quad |U_e^-| = k_e^-,\quad d^+_u\ge 2 \, \text{ for all } \,  u \in U_e^+, \quad d^-_u\geq 2\, \text{ for all } \, u\in U_e^-.
\] 
Then $Q(\dvec^+,\kvec^+,\dvec^-,\kvec^-)$ can be estimated from above by summing $\B_0(\hat{\dvec}^+,\hat{\kvec}^+,\hat{\dvec}^-,\hat{\kvec}^-)$ over all such pairs $e,f$ and sets $U_e^+,U_e^-$, where $\hat{\dvec}^+,\hat{\kvec}^+,\hat{\dvec}^-,\hat{\kvec}^-$ are constructed from $\dvec^+,\kvec^+,\dvec^-,\kvec^-$ as follows:
\begin{itemize}
    \item 
The vector $\hat{\dvec}^+$ matches $\dvec^+$ in every component except those corresponding to $u\in U_e^+$, where $\hat{d}_u^+ = d_u^+-2$. 
\item Similarly the vector $\hat{\dvec}^-$ matches $\dvec^-$ in every component except those corresponding to  $u\in U_e^-$, where $\hat{d}_u^- = d_u^- - 2$.
\item Finally $\hat{\kvec}^+$ matches $\kvec^+$ (respectively, $\hat{\kvec}^-$ matches $\kvec^-$) in every component except those corresponding to $e,f$, where $\hat{k}_e^+ = \hat{k}_f^+ = \hat{k}_e^- = \hat{k}_f^- = 0$.
\end{itemize}
The proof has three steps.

\medskip

\noindent {\bf Step 1:}\  Compare $\B_0(\hat{\dvec}^+,\hat{\kvec}^+,\hat{\dvec}^-,\hat{\kvec}^-)$ with $\B(\hat{\dvec}^+,\hat{\kvec}^+,\hat{\dvec}^-,\hat{\kvec}^-)$. 

\smallskip
\noindent Let $\hat{M}^+, \hat{M}^-, \hat{\Delta}$ correspond to $M^+,M^-,\Delta$ in the natural way. Then $M^{+}-2\Delta\le \hat{M}^+ \le M^+$ and $M^{-}-2\Delta\le \hat{M}^- \le M^-$, implying that
\[
\hat{M}^+ = \bigl(1+O(\Delta/M^+)\bigr) M^+ \quad \text{ and }\quad M^-=\bigl(1+O(\Delta/M^-)\bigr)M^-.
\]
As the value of any component of the vectors was only decreased we also have $\hat{\Delta}\le \Delta$.
In addition, for $\hat{\kvec}=(\hat{\kvec}^+,\hat{\kvec}^-)$ we have
\[
R(\kvec)-2\Delta^2\le R(\hat{\kvec}) \le R(\kvec)
\]
as two terms of $R(\kvec)$ which each may have been as large as $\Delta^2$ have been replaced by 0 in $R(\hat{\kvec})$.
Similarly, for $\hat{\dvec}=(\hat{\dvec}^+,\hat{\dvec}^-)$ we have
\[
R(\dvec)- 4\Delta^2 \le R(\hat{\dvec})\le R(\dvec)
\]
as decreasing any $d_u^+$ or $d_u^-$ by 2 decreases $R(\dvec)$ by at most $2\Delta$, and at most $2\Delta$ values have been changed. Therefore, using (\ref{eq:Rbounds}),
\[
\frac{R(\hat{\dvec})R(\hat{\kvec})}{\hat{M}^+\hat{M}^-}=\frac{R(\dvec)R(\kvec)}{M^{+}M^{-}}+O\left(\frac{\Delta^2 \bigl(R(\dvec)+R(\kvec)\bigr)}{M^+M^-}\right)=\frac{R(\dvec)R(\kvec)}{M^{+}M^{-}}+O\left(\frac{\Delta^3 }{\max\{M^+,M^-\}}\right).
\]
Combining this with Corollary~\ref{cor:B0}, we deduce that 
\begin{align}
&\B_0(\hat{\dvec}^+,\hat{\kvec}^+,\hat{\dvec}^-,\hat{\kvec}^-)\notag \\
&=\B(\hat{\dvec}^+,\hat{\kvec}^+,\hat{\dvec}^-,\hat{\kvec}^-)\, \exp\left(-\frac{R(\dvec)R(\kvec)}{M^+M^-}+O\left(\frac{\Delta^4}{\max\{M^+,M^-\}}+\frac{\Delta^2}{\min\{M^+,M^-\}}\right)\right)\notag \\
&= \B(\hat{\dvec}^+,\hat{\kvec}^+,\hat{\dvec}^-,\hat{\kvec}^-)\, \exp\left(-\frac{R(\dvec)R(\kvec)}{M^+M^-} + o(1)\right). \label{eq:end-step-1}
\end{align}

\noindent {\bf Step 2:}\ Compare $B(\hat{\dvec}^+,\hat{\kvec}^+,\hat{\dvec}^-,\hat{\kvec}^-)$ with $B(\dvec^+,\kvec^+,\dvec^-,\kvec^-)$.

\smallskip
\noindent Note, $\B(\hat{\dvec}^+,\hat{\kvec}^+,\hat{\dvec}^-,\hat{\kvec}^-)$ equals the number of directed bipartite graphs $G$ in $\calB(\dvec^+,\kvec^+,\dvec^-,\kvec^-)$ which contain all arcs in 
\[ \{(u,e), (u,f) : u\in U_e^+\} \cup \{(e,u), (f,u): u \in U_e^-\}.
\]
This follows since $\hat{k}_e^+=\hat{k}_e^-=\hat{k}_f^+=\hat{k}_f^-=0$, so
none of the above arcs can be present in any element of $\calB(\hat{\dvec}^+,\hat{\kvec}^+,\hat{\dvec}^-,\hat{\kvec}^-)$.
Applying Theorem~\ref{thm:subgraph} to $G^+$ and $G^-$ gives, for $n$ sufficiently large,
\begin{align}
\frac{\B(\hat{\dvec}^+,\hat{\kvec}^+,\hat{\dvec}^-,\hat{\kvec}^-)}{\B(\dvec^+,\kvec^+,\dvec^-,\kvec^-)}&\le \frac{(k_e^+!)^2(k_e^-!)^2 \prod_{u\in U_e^+} (d_u^+)_2 \prod_{u\in U_e^-} (d_u^-)_2}{\big(M^+-O(\Delta^2)\big)_{2k_e^+}\, \big(M^- -O(\Delta^2\big)_{2k_e^-}} \notag\\
&\leq (k_e^+!)^2(k_e^-!)^2 \prod_{u\in U_e^+} (d_u^+)_2 \prod_{u\in U_e^-} (d_u^-)_2\, \left(\frac{2}{M^+}\right)^{2k_e^+}\, \left(\frac{2}{M^-}\right)^{2k_e^-}. 
\label{eq:B-ratios}
\end{align}

\medskip

\noindent {\bf Step 3:}\  Sum over all choices of $e$, $f$, $U_e^+$, $U_e^-$.

\smallskip
\noindent First we will sum over all choices of $U_e^+$ and $U_e^-$ of the appropriate size, ignoring the disjointness condition for an upper bound. 
Observe that $(d_u^+)_2\leq \Delta d_u^+$ and $(d_u^-)_2\leq \Delta d_u^-$ for all $u\in V$, and that
\[
 k_{e}^+! \sum_{\substack{U_e^+\subseteq V\\|U_e^+|=k_e^+}}  \prod_{u\in U_e^+} d_u^+ \leq 
\bigg(\sum_{u\in V} d_u^+\bigg)^{k_e^+} = \big(M^+\big)^{k_e^+}
 \,\,\, \mbox{and}\,\,\, k_e^-! \sum_{\substack{U_e^-\subseteq V\\|U_e^-|=k_e^-}} \prod_{u\in U_e^-} d_u^- \leq \big( M^-\big)^{k_e^-}.
 \]
Using these inequalities with (\ref{eq:B-ratios}) gives
\[ 
\sum_{\substack{U_e^+\subseteq V\\|U_e^+|=k_e^+}} 
\sum_{\substack{U_e^-\subseteq V\\|U_e^-|=k_e^-}} 
\frac{\B(\hat{\dvec}^+,\hat{\kvec}^+,\hat{\dvec}^-,\hat{\kvec}^-)}{\B(\dvec^+,\kvec^+,\dvec^-,\kvec^-)}
\leq (k_e^+!)(k_e^-!)
\, \left(\frac{4\Delta}{M^+}\right)^{k_e^+}\, \left(\frac{4\Delta}{M^-}\right)^{k_e^-}. 
\]
 Under our assumptions, the right hand side is maximised when $k_e^+$ and $k_e^-$ are as small as possible. Due to our condition that $\kappa\geq 3$, 
 for an upper bound we can assume that either $k_e^+=2$ and $k_{e}^-=1$ or $k_e^+=1$ and $k_e^-=2$. 
Hence we obtain
\begin{align}
\sum_{\substack{U_e^+\subseteq V\\|U_e^+|=k_e^+}} 
\sum_{\substack{U_e^-\subseteq V\\|U_e^-|=k_e^-}} 
\frac{\B(\hat{\dvec}^+,\hat{\kvec}^+,\hat{\dvec}^-,\hat{\kvec}^-)}{\B(\dvec^+,\kvec^+,\dvec^-,\kvec^-)}
= O\left(\frac{\Delta^3}{M^+\, M^-\, \min\{M^+,M^-\}}\right).
\label{eq:sum-over-U}
\end{align}

Finally, for a worst case we assume that all $\binom{|E|}{2}$ choices of unordered pairs of hyperarc
$\{ e,f\}$ have the same head and tail sizes as each other. Since the tail and head of
each hyperarc is nonempty, it follows that $|E|\leq \min\{M^+,M^-\}$.  Hence, summing the above expression over all pairs $\{e,f\}$ and combining with Corollary~\ref{cor:B0}, (\ref{eq:end-step-1}) and (\ref{eq:sum-over-U}) gives
\begin{align*}
    & Q(\dvec^+,\kvec^+,\dvec^-,\kvec^-) \\
    &\leq 
    \sum_{\substack{e,f\in E\\ e\neq f}} \,\sum_{\substack{U_e^+, U_e^- \subseteq V\\|U_e^+|=k_e^+, |U_e^-|=k_e^-}}\, \B_0(\hat{\dvec}^+,\hat{\kvec}^+,\hat{\dvec}^-,\hat{\kvec}^- )\\
    &= \exp\left( - \frac{R(\dvec)R(\kvec)}{M^+\, M^-} + o(1)\right) \, 
     \sum_{\substack{e,f\in E\\ e\neq f}} \,\sum_{\substack{U_e^+, U_e^- \subseteq V\\|U_e^+|=k_e^+, |U_e^-|=k_e^-}}\, B(\hat{\dvec}^+,\hat{\kvec}^+,\hat{\dvec}^-,\hat{\kvec}^-)\\
        &= \exp\left( - \frac{R(\dvec)R(\kvec)}{M^+\, M^-}\right) \, B(\dvec^+,\kvec^+,\dvec^-,\kvec^-)\,
        \sum_{\substack{e,f\in E\\ e\neq f}}
    O\left(\frac{\Delta^3}{M^+\, M^-\, \min\{ M^+,\, M^-\}}\right)
    \\
    &= \B_0(\dvec^+,\kvec^+,\dvec^-,\kvec^-)\, \sum_{\substack{e,f\in E\\ e\neq f}} O\left(\frac{\Delta^3}{M^+\, M^-\, \min\{ M^+,\, M^-\}}\right)\\
    &=
    \B_0(\dvec^+,\kvec^+,\dvec^-,\kvec^-)\, O\left(\frac{\Delta^3}{\max\{M^+,M^-\}}\right),
\end{align*}
completing the proof of Lemma~\ref{lem:nomults}. 
\end{proof}

We can now put all the pieces together to prove our main result.

\begin{proof}[Proof of Theorem~\ref{thm:main}]
First note that 
\[ P(\dvec^+,\kvec^+,\dvec^-,\kvec^-) = 1 - \frac{Q(\dvec^+,\kvec^+,\dvec^-,\kvec^-)}{\B_0(\dvec^+,\kvec^+,\dvec^-,\kvec^-)} = 1 - O\left(\frac{\Delta^3}{\max\{M^+,M^-\}}\right)
\]
using Lemma~\ref{lem:nomults}.
Hence, combining (\ref{eq:L-B}), Corollary~\ref{cor:B0} and Lemma~\ref{lem:nomults} shows that
\begin{align*}
&L(\dvec^+,\kvec^+,\dvec^-,\kvec^-)\\
&= \exp\left( - \frac{R(\dvec)R(\kvec)}{M^+M^-} + O\left(\frac{\Delta^2}{\min\{M^+,M^-\}} + \frac{\Delta^4}{\max\{M^+,M^-\}}\right)\right)\, \B(\dvec^+,\kvec^+,\dvec^-,\kvec^-).
\end{align*}
The proof is completed using (\ref{eq:L-H-conversion}) and Corollary~\ref{cor:size-B}.
\end{proof}


\end{document}